\documentclass[11pt,a4paper]{article}
\usepackage[margin=25mm]{geometry}
\usepackage[T1]{fontenc}
\usepackage{lmodern,microtype,amsmath,amssymb,amsthm,booktabs}
\usepackage[hidelinks]{hyperref}
\hypersetup{
  pdftitle={Cyclotomic numerical semigroups with four minimal generators},
  pdfauthor={Yexu Xing}
}
\newtheorem{theorem}{Theorem}[section]
\newtheorem{lemma}[theorem]{Lemma}
\newtheorem{corollary}[theorem]{Corollary}
\newcommand{\N}{\mathbb Z_{\ge 0}}
\newcommand{\Q}{\mathbb Q}
\newcommand{\D}{\Delta}
\newcommand{\coeff}[1]{[t^{#1}]}
\title{Cyclotomic numerical semigroups\\[2pt]with four minimal generators}
\author{Yexu Xing}
\date{}

\begin{document}
\maketitle

\begin{abstract}
A numerical semigroup with four minimal generators is cyclotomic if and only
if its defining ideal is generated by three binomials. We prove the forward
implication by computing one coefficient of the formal logarithm of the
graded Hilbert numerator. For a symmetric numerical semigroup with four minimal generators and five minimal
relations, let \(F\) be its Frobenius number and let \(\Delta\) be \(F\) plus
the sum of its generators. The coefficient at \(\Delta\) is exactly \(4\).
The reciprocal roots of the numerator therefore have \(\Delta\)-th power
sum \(-4\Delta\), although there are only \(\Delta\) roots. The coefficient
calculation uses the eight nonzero entries above the diagonal in the
middle differential of the Bresinsky resolution and the gap at \(F\).
\end{abstract}

\section{The four-generator question}\label{sec:main}
Let \(S=\langle n_1,n_2,n_3,n_4\rangle\subset\N\) be minimally generated
by four positive integers of greatest common divisor one. Set
\(H_S(t)=\sum_{s\in S}t^s\), and write \(F=\max(\N\setminus S)\).
The semigroup polynomial and its values at the endpoints are given by
\[
 P_S(t)=(1-t)H_S(t)=1-(1-t)\sum_{g\in\N\setminus S}t^g,
 \qquad P_S(0)=P_S(1)=1 .
\]
In particular \(P_S\) is monic of degree \(F+1\). We call \(S\)
\emph{cyclotomic} if every complex root of \(P_S\) is a root of unity.
Let \(R=\Q[x_1,x_2,x_3,x_4]\) with \(\deg x_i=n_i\), and denote by \(I_S\)
the kernel of \(R\to\Q[t]\), \(x_i\mapsto t^{n_i}\).

\begin{theorem}\label{thm:main}
The semigroup \(S\) is cyclotomic if and only if \(I_S\) can be generated
by three binomials.
\end{theorem}

The image \(\Q[S]\) is integral over \(\Q[t^{n_1}]\), so \(I_S\) has height
three. The binomials equating monomials of the same weight span each
homogeneous part of \(I_S\). They therefore supply a minimal homogeneous
generating set, by taking a basis modulo \((x_1,x_2,x_3,x_4)I_S\).
Thus the three-binomial condition is equivalent to \(\Q[S]\) being a
complete intersection.

Ciolan, Garc\'ia-S\'anchez and Moree conjectured this equivalence for
every embedding dimension and established it through dimension three
\cite[Conjecture~1 and Lemma~7]{CGM}. Moscariello and Sammartano record
the four-generator case in \cite[\S7, Problem~27]{MS}. The Hilbert
numerator already appears in \cite[\S5]{CGM}. Zhang uses a reciprocal-root
power sum of \(P_S\) to treat a family defined by Rosales
\cite[Lemma~2.2 and Proposition~3.1]{Zhang}. Here the decisive index is
the last shift of the four-generator resolution. Its coefficient is the
same for every symmetric four-generator semigroup that requires five
relations.

\section{Degrees in the five-relation resolution}\label{sec:shifts}
Put \(\D=F+\sum_{i=1}^4 n_i\) and define the Hilbert numerator
\begin{equation}\label{eq:K}
 K(t)=H_S(t)\prod_{i=1}^4(1-t^{n_i})
     =P_S(t)\,\frac{\prod_{i=1}^4(1-t^{n_i})}{1-t}.
\end{equation}
It has constant term \(1\), degree \(\D\), and leading coefficient \(-1\).
The quotient multiplying \(P_S\) is a polynomial all of whose roots are
roots of unity.

Assume throughout this section and Section~\ref{sec:gap} that \(S\) is
symmetric and \(\Q[S]\) is not a complete intersection. Symmetry means
that exactly one of \(u,F-u\) belongs to \(S\) for every integer \(u\).
By Bresinsky's classification \cite{Bresinsky}, in the notation of
\cite[Theorem~3]{BFS}, the generators can be ordered so that the five
minimal relations have the following two monomials. This ordering is
cyclic and need not be increasing.
\begin{equation}\label{eq:rels}
\begin{array}{c@{\qquad}ll}
\toprule
i&\multicolumn{2}{c}{\text{Monomials of relation }i}\\
\midrule
1& x_1^{\alpha_1} & x_3^{\alpha_{13}}x_4^{\alpha_{14}}\\
2& x_2^{\alpha_2} & x_1^{\alpha_{21}}x_4^{\alpha_{24}}\\
3& x_3^{\alpha_3} & x_1^{\alpha_{31}}x_2^{\alpha_{32}}\\
4& x_4^{\alpha_4} & x_2^{\alpha_{42}}x_3^{\alpha_{43}}\\
5& x_2^{\alpha_{32}}x_4^{\alpha_{14}} &
   x_1^{\alpha_{21}}x_3^{\alpha_{43}}\\
\bottomrule
\end{array}
\end{equation}
All eight mixed exponents are positive, and
\begin{equation}\label{eq:alphas}
\alpha_1=\alpha_{21}+\alpha_{31},\quad
\alpha_2=\alpha_{32}+\alpha_{42},\quad
\alpha_3=\alpha_{13}+\alpha_{43},\quad
\alpha_4=\alpha_{14}+\alpha_{24}.
\end{equation}
Write \(d_i\) for the common weight in row \(i\) and \(b_i=\D-d_i\).
The minimal resolution gives
\begin{equation}\label{eq:resolution}
 K(t)=1-\sum_{i=1}^5 t^{d_i}+\sum_{i=1}^5 t^{b_i}-t^\D .
\end{equation}
In the alternating middle differential, an entry in position \(i,j\)
has weight \(b_j-d_i=\D-d_i-d_j\) whenever it is nonzero. In particular,
writing \(r_{ij}=\D-d_i-d_j\), its eight nonzero upper-triangular
entries give the following weights
\cite[Theorem~4 and Corollaries~13--14]{BFS}.
\begin{equation}\label{eq:residual}
\begin{array}{c|cccc|cccc}
 ij&12&23&34&14&15&25&35&45\\ \hline
 r_{ij}&\alpha_{43}n_3&\alpha_{14}n_4&\alpha_{21}n_1&\alpha_{32}n_2&
 \alpha_{24}n_4&\alpha_{31}n_1&\alpha_{42}n_2&\alpha_{13}n_3
\end{array}
\end{equation}
Each \(b_i\) is consequently the sum of some \(d_j\) and a positive
element of \(S\). In fact there are at least three distinct such indices
\(j\) for each \(i\). In particular all ten interior shifts are positive
and strictly smaller than \(\D\).

Each weight in \eqref{eq:residual} has a unique factorization in the
four generators.
Indeed, setting all variables other than \(x_j\) to zero sends \(I_S\)
to \((x_j^{\alpha_j})\) by \eqref{eq:rels}. If \(h n_j\), with
\(0<h<\alpha_j\), had another factorization, cancellation of all
common copies of \(n_j\) would give a binomial equating a positive
power of \(x_j\) below \(\alpha_j\) to a monomial in the other
variables. Its image contradicts the displayed ideal. Equation
\eqref{eq:alphas} gives \(0<\alpha_{ij}<\alpha_j\) for each multiplier
of \(n_j\) in \eqref{eq:residual}.

\section{Representing the top degree}\label{sec:gap}
A monomial of weight \(\D\) cannot involve all four variables. Dividing
such a monomial by \(x_1x_2x_3x_4\) would give one of weight
\(\D-\sum n_i=F\), which is a gap. This observation restricts sums of
relation degrees because each summand admits both factorizations in
\eqref{eq:rels}.

\begin{lemma}\label{lem:ordinary}
The top degree \(\D\) does not belong to \(\langle d_1,\ldots,d_5\rangle\).
\end{lemma}
\begin{proof}
Suppose that \(\D\) is a sum of relation degrees, with repetition
allowed. First consider two distinct types \(d_i,d_j\) in this sum.
For \(\{i,j\}=\{1,3\}\) or \(\{2,4\}\), their mixed monomials in
\eqref{eq:rels} together involve all four variables, a contradiction.
For every other pair, the remaining summands must give \(r_{ij}\).
This is a positive weight with a unique factorization by
\eqref{eq:residual}. The remaining sum must be nonempty. Each of
its summands has two different factorizations, and adding any fixed
factorization of the other summands preserves that difference. This
contradicts uniqueness.

It remains to rule out \(\D=k d_i\). Since \(d_i<\D\), we have \(k\ge2\).
For \(i=5\), one copy of each monomial in row 5 already uses all four
variables. For the first four rows a cyclic permutation of \(1,2,3,4\)
preserves the relations up to reordering and a change of sign in the
fifth. We may therefore take \(i=1\). Put \(a=\alpha_{13}\) and
\(b=\alpha_{14}\), so \(k d_1=ka n_3+kb n_4\).

If \(ka>\alpha_3\), replace \(\alpha_3 n_3\) by
\(\alpha_{31}n_1+\alpha_{32}n_2\). Positive copies of both \(n_3\)
and \(n_4\) remain. If \(kb>\alpha_4\), replace \(\alpha_4 n_4\) by
\(\alpha_{42}n_2+\alpha_{43}n_3\). We can now also replace
\(\alpha_3 n_3\) by \(\alpha_{31}n_1+\alpha_{32}n_2\), since the
number of copies of \(n_3\) left after this substitution is
\[
 ka+\alpha_{43}-\alpha_3=(k-1)a>0 .
\]
There are still \(kb-\alpha_4>0\) copies of \(n_4\). Both cases yield
a factorization of \(\D\) involving all four generators. In the remaining
case \(ka\le\alpha_3\) and \(kb\le\alpha_4\). Then
\[
 k d_1\le d_3+d_4=\D-\alpha_{21}n_1<\D
\]
by \eqref{eq:residual}. This finishes the proof.
\end{proof}

The first shifts \(d_i\) form an antichain in the semigroup order,
\begin{equation}\label{eq:antichain}
 d_i-d_j\notin S\qquad(i\ne j).
\end{equation}
This follows from \cite[Lemma~19 and Theorem~27]{BFS} and can also be
read directly from \eqref{eq:rels}. At a fixed weight, take the monomials
of that weight as vertices. Join two when their difference is a monomial
multiple of a listed relation, up to sign. Since these relations generate
\(I_S\), the edge differences span all differences in the fiber, so the
graph is connected. At weight \(d_i\), neither term of row \(i\)
is divisible by a term of another row. The six mixed supports are
distinct pairs, and each mixed exponent is below the corresponding
pure exponent by \eqref{eq:alphas}. The two terms of row \(i\) thus
form the entire fiber and are coprime. If \(d_i=d_j\) for \(j\ne i\),
row \(j\) would give additional monomials in this fiber. If \(d_i-d_j\)
were positive and in \(S\), multiplying both terms of row \(j\) by a
monomial of weight \(d_i-d_j\) would give two distinct monomials in the
fiber with a common factor. This proves \eqref{eq:antichain}.

\begin{lemma}\label{lem:dual}
The only multisets of interior shifts \(d_j,b_j\) with sum \(\D\) are
the five pairs \(\{d_i,b_i\}\).
\end{lemma}
\begin{proof}
Lemma~\ref{lem:ordinary} forces an occurrence of some \(b_i\), leaving
a sum of weight \(d_i\). If the remainder contained a \(b_j\), choose
\(k\ne i\) such that \(b_j-d_k\in S\) using \eqref{eq:residual}.
Every summand belongs to \(S\), so \(d_i-d_k\in S\), contrary to
\eqref{eq:antichain}. For the same reason a first shift \(d_j\) in the
remainder must have \(j=i\). Since \(d_i>0\), exactly one copy occurs.
Conversely, each of the five indicated pairs has weight \(\D\).
\end{proof}

\section{The reciprocal-root power sum}\label{sec:power}
For a polynomial \(A\) of positive degree with \(A(0)=1\), let
\(s_N(A)=\sum_{\rho}\rho^{-N}\) sum over its complex roots with
multiplicity. Expanding \(\log A(t)=\sum_\rho\log(1-t/\rho)\) at zero
gives \(s_N(A)=-N\coeff{N}\log A(t)\) for every \(N\ge1\). If all
roots lie on the unit circle, then \(\lvert s_N(A)\rvert\le\deg A\).

\begin{theorem}\label{thm:moment}
For every symmetric numerical semigroup with four minimal generators
whose defining ideal requires five binomials,
\[
 \coeff{\D}\log K(t)=4,\qquad s_\D(K)=-4\D .
\]
In particular \(P_S\) has a root off the unit circle.
\end{theorem}
\begin{proof}
Write \(Z=K-1\) in \eqref{eq:resolution} and expand
\(\log K=Z-Z^2/2+Z^3/3-\cdots\) as a formal power series.
The linear term contributes \(-1\) at weight \(\D\). By
Lemma~\ref{lem:dual}, every nonlinear contribution uses one \(d_i\)
and its complementary \(b_i\). For each \(i\) the two orders in
\(Z^2\) both have coefficient \(-1\). Thus \(-Z^2/2\) contributes
\(5\), and all higher powers contribute zero. This proves both
identities. The absolute power sum \(4\D\) exceeds \(\deg K=\D\).
Equation~\eqref{eq:K} transfers the resulting off-circle root to \(P_S\).
\end{proof}

\begin{proof}[Proof of Theorem~\ref{thm:main}]
Every cyclotomic numerical semigroup is symmetric
\cite[Theorem~1]{CGM}. If \(I_S\) were not a complete intersection,
Bresinsky's classification would give the five relations in
\eqref{eq:rels}, contradicting Theorem~\ref{thm:moment}.
Conversely, three homogeneous binomial generators form a regular
sequence, since their number equals the height in the polynomial ring.
For their degrees \(e_1,e_2,e_3\), the Koszul resolution gives
\(K(t)=\prod_{j=1}^3(1-t^{e_j})\). Every root of \(P_S\) is a root of
this product by \eqref{eq:K}, hence is a root of unity.
\end{proof}

The identity can also be read directly on the semigroup polynomial.
Taking the logarithm in \eqref{eq:K} under the hypotheses of
Theorem~\ref{thm:moment} gives
\[
 s_\D(P_S)=1-4\D-\sum_{i:\,n_i\mid \D}n_i .
\]
For an example at multiplicity eight, take \(S=\langle8,9,10,11\rangle\).
The least elements of \(S\) in the residue classes \(0,\ldots,7\) modulo
eight are \(0,9,10,11,20,21,22,31\), respectively. Each element of
\(S\) is uniquely one of these representatives plus a nonnegative
multiple of eight. Thus \(F=31-8=23\), \(\D=61\), and
\[
 K(t)=(1+t^9+t^{10}+t^{11}+t^{20}+t^{21}+t^{22}+t^{31})
       (1-t^9)(1-t^{10})(1-t^{11}).
\]
Multiplying gives the five negative interior exponents \(18,19,20,32,33\)
and the five positive ones \(61-d_i\). No sum of two negative interior
exponents is \(61\). The sum of three is at most \(60\) if all are
at most \(20\), and at least \(32+2\cdot18>61\) otherwise. Four already
sum to at least \(72\). This illustrates Lemma~\ref{lem:ordinary}
even though two relation degrees exceed \(\D/2\). The complementary
pairs give \(s_{61}(K)=-244\), and the formula for \(P_S\) gives
\(s_{61}(P_S)=-243\).

\section{Multiplicity and higher embedding dimension}\label{sec:further}
\begin{corollary}\label{cor:mult}
The multiplicities of cyclotomic numerical semigroups with four
minimal generators are exactly the integers \(m\ge8\).
\end{corollary}
\begin{proof}
The multiplicity is the least positive element of \(S\). Li and Zhang
\cite[Theorem~1.4]{LZ} show that for every \(e\ge2\), a complete
intersection of embedding dimension \(e\) and multiplicity \(m\)
exists exactly when \(m\ge2^{e-1}\). Apply this statement with \(e=4\)
and use both directions of Theorem~\ref{thm:main}.
\end{proof}

Moscariello and Sammartano \cite[\S7, after Problem~29]{MS} trace the
complete-intersection realization theorem to Watanabe
\cite[Theorem~1 and the subsequent remark]{Watanabe}.
At the threshold, \(\langle8,12,18,27\rangle\) gives a cyclotomic
example. The three relations equate \(x_2^2,x_3^2,x_4^2\) with
\(x_1^3,x_2^3,x_3^3\), respectively. Reduction in descending
variable order gives representatives \(x_1^q x_2^u x_3^v x_4^w\)
with \(u,v,w\in\{0,1\}\). Their tails have distinct residues
\(4u+2v+3w\bmod 8\). Parity determines \(w\), reduction modulo four
then determines \(v\), and the remaining residue determines \(u\).
Thus equal weights give the same representative, proving that the
three binomials generate the ideal.

Theorem~\ref{thm:main} also settles the conjecture for any five-generator
semigroup that is a proper gluing. Indeed, let \(S=aS_1+bS_2\) be a gluing
in which the scaled minimal generating sets of \(S_1\) and \(S_2\)
partition the five minimal generators of \(S\). The gluing formula
\cite[\S5, (20)]{CGM} gives
\[
 P_S(t)=P_{\langle a,b\rangle}(t)\,P_{S_1}(t^a)\,P_{S_2}(t^b).
\]
If \(S\) is cyclotomic, every root of either \(P_{S_i}\) is a power
of a root of \(P_S\), hence a root of unity. The two factors have at
most four minimal generators each, so Theorem~\ref{thm:main} and the
lower-dimensional cases \cite[Lemma~7]{CGM} make both complete
intersections. Their gluing is a complete intersection by Delorme's
theorem \cite[\S5]{CGM}.

For a monic polynomial \(A\in\mathbb Z[t]\) with \(A(0)=1\), a root
that is not a root of unity can always be detected by a power sum.
Indeed, Kronecker's theorem \cite[Lemma~1]{CGM} then gives a root of
modulus greater than one. The product of all root moduli is one, so
there is also a root inside the unit disk. At this root the function
\[
 -\frac{tA'(t)}{A(t)}=\sum_{N\ge1}s_N(A)t^N
\]
has a pole. Bounded coefficients would make the series converge throughout
the open unit disk. Thus \(\lvert s_N(A)\rvert>\deg A\) for some \(N\).
Conversely, if every root is a root of unity, the triangle inequality
gives \(\lvert s_N(A)\rvert\le\deg A\) for all \(N\).

For five minimal generators \(n_1,\ldots,n_5\), put
\(\D=F+\sum_i n_i\) and \(K(t)=H_S(t)\prod_i(1-t^{n_i})\).
As in \eqref{eq:K}, \(K\) has degree \(\D\) and the same roots off
the unit circle as \(P_S\). Its leading coefficient is now \(1\).
The gap at \(F\) still excludes a full-support factorization of \(\D\).
Nevertheless, the top-degree power sum can vanish.

Take \(S=\langle10,15,16,19,22\rangle\). Its least representatives
modulo \(10\), in residue order, are
\(W=(0,31,22,53,34,15,16,37,38,19)\). Hence \(F=43\) and there are
\(\sum_{w\in W}\lfloor w/10\rfloor=22\) gaps, so \(S\) is symmetric.
Its multiplicity \(10<2^4\) excludes a complete intersection
\cite[Theorem~1.4]{LZ}. Here \(\D=125\) and
\[
 K(t)=\left(\sum_{w\in W}t^w\right)
      (1-t^{15})(1-t^{16})(1-t^{19})(1-t^{22}).
\]
Put \(Z=K-1\). Its lowest term is \(-t^{30}\). At degree \(125\) the
coefficients of \(Z,Z^2,Z^3,Z^4\) are \(1,22,33,4\), respectively,
so \(\coeff{125}\log K=1-22/2+33/3-4/4=0\). At degree \(73\) there
is no linear term and the only contributing pairs are \(32+41\) and
\(35+38\), each formed from two coefficients equal to \(-1\). Thus
\(\coeff{73}\log K=-2\), giving \(s_{125}(K)=0\) and \(s_{73}(K)=146\).

Among the symmetric semigroups with five minimal generators,
multiplicity at most \(15\), and largest generator at most \(32\),
there are \(746\) examples. All have multiplicity below the
complete-intersection threshold. Enumerating the increasing generator
tuples, computing least representatives by shortest paths modulo the
multiplicity, and extracting the coefficients of \(-tK'/K\) gives
\(\lvert s_\D(K)\rvert>\D\) in \(690\) cases. Each of the other \(56\)
has \(\lvert s_N(K)\rvert>\D\) for some \(N<\D\). A reader of an earlier
version suggested the upper bound in the following question.

\smallskip
\noindent\textit{Question.}
If \(S\) is a symmetric numerical semigroup with five minimal
generators and \(\Q[S]\) is not a complete intersection,
must there be an integer \(1\le N\le\D\) with \(\lvert s_N(K)\rvert>\D\)?
An affirmative answer would prove the five-generator case of the
conjecture, since every cyclotomic numerical semigroup is symmetric.

\section*{Use of AI}
OpenAI Codex assisted in searching the literature, developing the proof,
drafting and editing the manuscript, and writing Lean code.
The author takes responsibility for the manuscript.

\bigskip
\noindent\begin{minipage}{\textwidth}
\small
\textsc{Fudan University, Shanghai, China}\\
\textsc{Shanghai Institute for Mathematics and Interdisciplinary Sciences (SIMIS),}\\
\textsc{Shanghai, China}\\[3pt]
\textit{Email address:}
\href{mailto:25114020022@m.fudan.edu.cn}{\nolinkurl{25114020022@m.fudan.edu.cn}}
\end{minipage}
\end{document}